\documentclass[11pt]{article}
\usepackage{amsbsy,amsthm,amsfonts,amssymb}
\usepackage{amsmath}
\usepackage{amsfonts}
\usepackage{verbatim}
\usepackage[mathscr]{eucal}
\newtheorem{theo}{Theorem}[section]
\newtheorem{lemma}[theo]{Lemma}
\newtheorem{prop}[theo]{Proposition}
\newtheorem{example}[theo]{Example}
\newtheorem{defi}[theo]{Definition}
\newtheorem{coll}[theo]{Corollary}

\newtheorem*{thm1}{Theorem~\ref{maintheo3}}
\newtheorem*{thm2}{Theorem~\ref{maintheo2}}
\newtheorem*{thm3}{Theorem~\ref{maintheo4}}
\newtheorem*{cor}{Corollary~\ref{ci}}

\newcommand{\mb}{\mathcal{B}}
\newcommand{\be}{\begin{equation}}
\newcommand{\ee}{\end{equation}}
\newcommand{\ba}{\begin{array}}
\newcommand{\ea}{\end{array}}

\begin{document}

\title{Some arithmetic properties of matroidal ideals}

\author{ Hung-Jen Chiang-Hsieh
\thanks{e-mail: hchiang@math.ccu.edu.tw} ~~and~~
 Hsin-Ju Wang \thanks{e-mail:  hjwang@math.ccu.edu.tw} \\
Department of Mathematics, National Chung Cheng University,\\
          Chiayi 621, Taiwan}

\date{}

\maketitle

\noindent{\bf Abstract.} In this paper, we study various
properties of matroidal ideals.

\section{Introduction}
Let $K$ be a field and $R=K[x_1, \dots, x_n]$ be the
polynomial ring in $n$ variables over $K$ with each $\deg~x_i=1$. %
If $u=x_1^{a_1}\cdots x_n^{a_n}$ is a monomial of $R$, then we
denote the support of $u$ by $supp(u)=\{\,x_i~|~a_i\neq 0\}$. For a
monomial ideal $I\subseteq R$, $G(I)$ is denoted for the set of its
unique minimal monomial generators. We call a monomial ideal $I$ a
\emph{matroidal ideal} if each member of $G(I)$ is square-free
(i.e., $I$ is reduced) and that the following exchange condition is
satisfied: for any $u=x_1^{a_1}\cdots x_n^{a_n}, v=x_1^{b_1}\cdots
x_n^{b_n}\in G(I)$, if $a_i>b_i$ for some $i$, then there exists
some $j$ with $a_j<b_j$
such that $x_ju/x_i \in G(I)$%
\footnote{$I$ is called a \emph{polymatroidal ideal} when not
requiring the square-free assumption(see\cite{hh}).  }

In other words, the set $\mb(I) =\{\,supp(u)\,|\,u\in G(I)\,\}$
satisfies the following exchange condition:
\begin{list}{(B)}
\item  If $B_1$ and $B_2$ are elements of $\mb(I)$ and $x \in B_1 -
B_2$, then there is an element $y\in B_2 - B_1$ such that
$(B_1-\{x\})\cup \{y\} \in \mb(I)$.
\end{list}

\medskip
\noindent It follows form \cite[Theorem~1.2.3]{ox} that there is a
``matroid" having $\mb(I)$ as its collection of bases (maximal
independent sets). Since each maximal independent set of a matroid
has the same cardinality (see \cite[Lemma~1.2.4]{ox}), each monomial
$u\in G(I)$ must be of the same degree, say $d$, and we call this
number $d$ the degree of the matroidal ideal $I$. 

The matroid theory is one of the most fascinating research area in
combinatorics which has many links to graphs, lattices, codes, and
projective geometry. For the interested reader, we refer the
textbooks \cite{ox} or \cite {we}. In this paper, we focus on some
arithmetic properties held by a matroidal ideal.

It is known that a matroidal ideal has linear quotients (cf.
\cite[Theorem~5.2]{ch}). We first discuss the linear quotient index
$q(I)$ of a matroidal ideal in section two and get the following
result.

\begin{thm1} Let $I$ be a matroidal ideal of degree $d$ in the
polynomial ring $R=K[x_1, \dots, x_n]$ with
$supp(I)=\{x_1,\dots,x_n\}$. Then $q(I)=n-d$.
\end{thm1}

\noindent With this result and the fact \cite[Corollary~1.6]{ht} we
obtain that the projective dimension of a matroidal ideal is
$pd_R(I)=n-d$.

An ideal is \emph{unmixed} if all its prime divisors are of the same
height. It is known that the Cohen-Macaulay ideals hold this
property. In section three, we discuss the unmixed matroid ideal and
find the relation between the height and the degree of a matroid
ideal as follows. (for the ``$*$" product of ideals, please see
Definition 2.1)

\begin{thm2}
\label{maintheo2} Let $I \subseteq K[x_1, \dots, x_n]$ be an unmixed
matroidal ideal of degree $d$ 
with $supp(I)=\{x_1,\dots,x_n\}$ and $n\geq 2$; then $h+d-1\leq
n\leq hd$, where $h$ is the height of $I$. In particular, $n=h+d-1$
if and only if $I$ is square-free Veronese; and $n=hd$ if and only
if $I=J_1*J_2* \cdots *J_d$, where each $J_i$ is generated by $h$
distinct variables.
\end{thm2}

For an ideal $I$, the minimal number of elements which generate $I$
up to radical is called the \emph{arithmetical rank} of this ideal
and is denoted by $ara\, I$. When this numerical invariant equals to
the height of $I$, we say that $I$ is a set-theoretic complete
intersection. We discuss the relation between the arithmetical rank
$ara\,I$ and the linear quotient index $q(I)$ of a matroidal ideal
in the final section. The main result we obtain is as below.

\begin{thm3}\label{maintheo4} Let $I$ be a matroidal ideal of degree $d$ of
a polynomial ring $R=K[x_1, \dots, x_n]$ and supp$(I)=\{x_1, \dots,
x_n \}$. Then $ara~I =q(I)+1$ if one of the following conditions
holds:
\begin{description}
\item{(i)} $I$ is square-free  Veronese;  \item{(ii)} $I=J_1 J_2 \cdots  J_d$
such that each $J_i$ is generated by $h$ distinct variables;
\item{(iii)} $d=2$,
\end{description} where $h$ is the height of $I$.
\end{thm3}

\noindent As a consequence of the above theorem, we have the
corollary.

\begin{cor} Let $I\subseteq K[x_1,\dots,x_n]$ be a matroidal ideal
such that supp$(I)=\{x_1, \dots, x_n \}$. Then $I$ is Cohen-Macaulay
if and only if it is a set-theoretic complete intersection.
\end{cor}

\section{Linear quotients and matroidal ideals}
Throughout, $R=K[x_1, \dots, x_n]$ is the polynomial ring in $n$
variables over a field $K$. By Cohen-Macaulay for an ideal $I$, we
mean that the quotient ring $R/I$ is Cohen-Macaulay. For a monomial
ideal $I$, we define the support of $I$ to be the set
$supp(I)={\bigcup}_{u\in G(I)} supp(u)$. In this section, we discuss
the linear quotient index $q(I)$ of a matroidal ideal. We first
recall the following definition.


\begin{defi}
We say that a monomial ideal $I\subseteq R$ has linear quotients if
there is an ordering $u_1, \dots, u_s$ of the monomials belonging to
$G(I)$ with deg$~u_1\leq  $deg$~u_2\leq \cdots \leq $deg$~u_s$ such
that, for each $2\leq j\leq s$, the colon ideal $\langle u_1, u_2,
\dots, u_{j-1}\rangle:u_j$ is generated by a subset of $\{x_1,
\dots, x_n\}$. \end{defi}

Let $I$ be a monomial ideal with linear quotients with respect to
the ordering $\{u_1, \dots, u_s\}$ of the monomials belonging to
$G(I)$. We write $q_j (I)$ for the number of variables which is
required to generate the colon ideal $\langle u_1, u_2, \dots,
u_{j-1}\rangle : u_j$. Let $q(I)=max\{ q_j(I)~|~2\leq j\leq s\}$.
From the fact \cite[Corollary~1.6]{ht} that the length of the
minimal free resolution of $R/I$ over $R$ is equal to $q(I)+1$, we
see that the index $q(I)$ is independent of the particular choice of
the ordering of the monomials which gives linear quotients.
Moreover, by the Auslander-Buchsbaum formula, we have
$depth\,R/I=n-q(I)-1$. It then follows from the equality $\dim
R/I=n-ht(I)$ that a monomial ideal $I$ with linear quotients
satisfies $ht(I)\leq q(I)+1$ and is Cohen-Macaulay if and only if
$ht(I)=q(I)+1$. We summarize the above as the following proposition.

\begin{prop}\label{qI}
Let $I$ be a monomial ideal of $R$ with linear quotients. Then
$ht(I)\leq q(I)+1$; and $I$ is Cohen-Macaulay if and only if
$\,ht(I)=q(I)+1$.
\end{prop}

As stated in the introduction, it is known that the matroidal ideals
have linear quotients. Therefore all the above discussion applied to
matroidal ideals. Next, we introduce two lemmas which are useful
later. In the sequel, we say that $I$ is a matroidal ideal of
$K[x_1,\dots,x_n]$ if $supp(I)=\{x_1,\dots,x_n\}$.




\begin{lemma}
\label{keylem3} Let $I \subseteq K[x_1, \dots, x_n]$ be a matroidal
ideal of degree $d$; 
and let $x$ and $y$ be variables in $R$ such that $xy\nmid u$ for
any $u\in G(I)$. If $xf\in G(I)$ for some monomial $f$ of degree
$d-1$, then $yf\in G(I)$.
\end{lemma}

\begin{proof}
Write $f=x_1\cdots x_{d-1}$. The assertion is clear if $d=2$ so we
may assume that $d\geq 3$. Let $g=y_1 \cdots y_{d-1}$ be a monomial
in $R$ different from $f$ such that $yg\in G(I)$ and $|supp(f)\cap
supp(g)|$ is maximal. We may assume that $y_i=x_i$ for $i=1, \dots,
k$. Suppose that $k\leq d-2$. Then by the definition of matroidal
ideal there are integers $i, j\geq k+1$ such that
$\frac{yg}{y_j}x_i\in I$, which contradicts to the choice of $g$.
Therefore, $f=g$ and the assertion holds.
\end{proof}

\begin{lemma}
\label{keylem2} Let $I \subseteq K[x_1, \dots, x_n]$ be a matroidal
ideal of degree $d$. 
If there are $d+1$ distinct variables $\{y, y_1, \dots,
y_d\}\subseteq \{x_1, \dots, x_n\}$ such that $f=y_1\cdots y_d\in
I$, then there exists an integer $i$ such that $\frac{f}{y_i} y\in
I$.
\end{lemma}

\begin{proof}
The assertion is clear if $d$ is small. We may assume that $d\geq
3$. Let $g=z_1 \cdots z_d$ be a monomial in $I$ different from $f$
such that $y\in supp(g)$ and $|supp(f)\cap supp(g)|$ is maximal. We
may assume that $z_i=y_i$ for $i=1, \dots, k$ and $z_d=y$. Suppose
that $k\leq d-2$. Then by the definition of matroidal ideal there
are integers $i, j\geq k+1$ such that $\frac{g}{z_j}x_i\in I$, which
contradicts to the choice of $g$. Therefore, $k=d-1$ and the
assertion holds.
\end{proof}

\begin{theo}
\label{maintheo3} Let $I$ be a matroidal ideal of degree $d$ of the
polynomial ring $R=K[x_1, \dots, x_n]$ with
$supp(I)=\{x_1,\dots,x_n\}$. Then $q(I)=n-d$.
\end{theo}

\begin{proof}
Since $I$ has linear quotients, 
there is an ordering $u_1, \dots, u_s$ of the monomials belonging to
$G(I)$ such that, for each $2\leq j\leq s$, the colon ideal $\langle
u_1, u_2, \dots, u_{j-1}\rangle:u_j$ is generated by a subset of
$\{x_1, \dots, x_n\}$.\\ To show the assertion, it is enough to show
that \be \label{eq1} \langle u_1, u_2, \dots,
u_{j-1}\rangle:u_j\subseteq \{x_1, \dots, x_n\}-supp(u_j)\ee for
each $2\leq j\leq s$ and
$$\langle u_1, u_2, \dots, u_{s-1}\rangle:u_s=\{x_1, \dots, x_n\}-supp(u_s).$$
Write $u_j=x_{i_1}\cdots x_{i_d}$. If $x_{i_t}\in \langle u_1, u_2,
\dots, u_{j-1}\rangle:u_j$ for some $t\leq d$, then $u_j\in \langle
u_1, u_2, \dots, u_{j-1}\rangle$ as $\langle u_1, u_2, \dots,
u_{j-1}\rangle$ is a square-free monomial ideal, a contradiction.
Thus, (\ref{eq1}) holds. By (\ref{eq1}), to finish the proof, it
suffices to show that $y\in \langle u_1, u_2, \dots,
u_{s-1}\rangle:u_s$ if $y\notin supp(u_s)$. However, this follows by
Lemma~\ref{keylem2} with $u_s=y_1\cdots y_d$.
\end{proof}

\begin{coll} \label{pd}
Let $I$ be a matroidal ideal of degree $d$ of the
polynomial ring $R=K[x_1, \dots, x_n]$. 
Then the projective dimension of the
ideal $I$ over $R$ is $pd_R(I)=n-d$.
\end{coll}

\begin{proof}
Since the length of the minimal free resolution of $R/I$ over $R$ is
$q(I)+1$ (see \cite[Corollary~1.6]{ht}), we obtain that
$pd_R(I)=pd_R(R/I)-1=q(I)=n-d$
\end{proof}

\section{Unmixed matroidal ideals}
An ideal is \emph{unmixed} if all its prime divisors are of the same
height. This property is held by a Cohen-Macaulay ideal. In this
section, we give characterizations of an unmixed matroidal ideal in
terms of its height, degree, and the number of variables.

We first recall one special kind of matroidal ideals, the
square-free  Veronese ideals.
\begin{example}
The square-free  Veronese ideal of degree $d$ in the variables
$\{x_1, \dots, x_n \}$ is the ideal which is generated by all
square-free monomials in $\{x_1, \dots, x_n \}$ of degree $d$. It is
easy to see that the square-free Veronese ideals are matroidal and
unmixed. In particular, from \cite[Theorem~4.2]{hh} one see that the
square-free Veronese ideals are the only case for Cohen-Macaulay
matroidal ideals.
\end{example}

We now give a characterization of matroidal ideal of degree $2$.
\begin{theo}
\label{keyprop} Let $I\subseteq K[x_1, \dots, x_n]$ be a matroidal
ideal of degree $2$ with $supp(I)=\{x_1,\dots,x_n\}$.
Then there are subsets $S_1, \dots, S_m$ of $\{x_1, \dots, x_n\}$
such that the following hold:
\begin{description}
\item{(i)} $m\geq 2$ and $|S_i|\geq 1$ for each $i$;
\item{(ii)} $S_i\cap S_j=\emptyset$ if $i\neq j$ and $\bigcup_{i=1}^m S_i=\{x_1,
\dots,x_n\}$;
\item{(iii)} if $x\in S_i$, $y\in S_j$ for $i\neq j$, then $xy\in
G(I)$;
\item{(iv)} if $x, y\in S_i$ for some $i$, then $xy\notin G(I)$.
\end{description}
Moreover, let $P_i$ be the  the prime ideals generated by the set
$\{x_1,\dots,x_n\}-S_i$ for each $i$. Then $P_1\cap P_2\cap \cdots
\cap P_m$ gives the primary decomposition of $I$.
\end{theo}

\begin{proof}
Let $t-1=|\{x_i~|~i\neq 1, x_ix_1\notin G(I)\}|$; then $1\leq t\leq
n-1$. Without loss of generality, we may assume that $x_1x_i\notin
G(I)$ if $i=2, \dots, t$ and $x_1x_i\in G(I)$ if $i=t+1, \dots, n$.
We first show the following two statements:
\\ (a) $x_ix_j\in G(I)$ if $i\leq t$ and $j\geq t+1$; \\ (b) $x_ix_j\notin G(I)$ if $i, j\leq
t$.

To show (a) holds, suppose on the contrary that $x_ix_j\notin G(I)$
for some $i\leq t$ and $j\geq t+1$. Since $x_i\in supp(I)$, there is
a variable $x_k$ such that $x_ix_k\in G(I)$. Moreover, $x_1x_j,
x_ix_k\in G(I)$ and $I$ is matroidal imply that either $x_ix_1$ or
$x_ix_j$ is in $G(I)$, a contradiction. Thus (a) holds.  For (b),
suppose on the contrary that $x_ix_j\in G(I)$ for some $i,j\leq t$.
Since $x_1x_n, x_ix_j\in G(I)$, it follows from the exchange
property of matroidal ideals that either $x_1x_i$ or $x_1x_j$
belongs to $G(I)$, a contradiction. Thus (b) holds. \par Let
$S_1=\{x_1, \dots, x_t\}$. Observe that $\{x_ix_j~|~i\leq
t,~and~j\geq t+1\}$ is a subset of $G(I)$. If $G(I)=\{x_ix_j~|~i\leq
t,~and~j\geq t+1\}$ then by setting $S_2=\{x_{t+1}, \dots, x_n\}$
and we are done. Therefore, we may assume that there are $j,k\geq
t+1$ such that $x_jx_k\in G(I)$. Let $I'$ be the monomial ideal in
$K[x_{t+1}\dots,x_n]$ generated by the set $G(I)-\{x_ix_j~|~i\leq
t,~and~j\geq t+1\}$. Then $supp(I')\subseteq \{x_{t+1}, \dots,
x_n\}$. In fact, $supp(I')=\{x_{t+1}, \dots, x_n\}$. For if not,
then there is a variable $x_l$ with $l\geq t+1$ such that $x_l\notin
supp(I')$. Since $x_lx_1, x_jx_k\in G(I)$ and $I$ is matroidal,
either $x_lx_j$ or $x_lx_k$ is in $G(I)$. Therefore, either $x_lx_j$
or $x_lx_k$ is in $G(I')$, a contradiction. We note that $I'$ is a
matroidal ideal of degree $2$ of the polynomial ring $K[x_{t+1},
\dots, x_n]$. Thus, the assertion follows by induction.


 Let $P_i$ be the prime ideals generated by the set
$\{x_1,\dots,x_n\}-S_i$. By the properties of $S_i$, it is easy to
see that $P_i=I:y$ for every $y\in S_i$. Therefore each $P_i$ is an
associate prime ideal of $I$. Let $w\in P_1\cap\dots \cap P_m$; then
$w\cdot y\in I$ whence $y\in \bigcup_{i=1}^{m} S_i$. It follows that
$(I:w)\supseteq\langle x_1,\dots,x_n\rangle$. Since $I$ is reduced,
$I$ has no embedded prime ideals. Therefore $I:w=R$ and so that
$w\in I$. Hence, $P_1\cap\dots \cap P_m=I$; and this completes the
proof.
\end{proof}

From the above theorem we see that the $S_i\,'s$ 
are uniquely determined. Moreover, if $I$ is unmixed then we have
that $|S_i|=|S_j|=n-ht(I)$ for all $i,j$.
Therefore we have the following corollary. 

\begin{coll}\label{d=2}
Let $I\subseteq K[x_1, \dots, x_n]$ be an unmixed matroidal ideal of
degree $2$; then one has $\frac{n}{2}\leq ht(I)\leq n-1 $. In
particular, $ht(I)=n-1$ if and only if $I$ is square-free Veronese;
and $ht(I)=\frac{n}{2}$ if and and only if $n$ is even and
$I=I_1*I_2$ such that each $I_i$ is generated by $\frac{n}{2}$
distinct variables.
\end{coll}

\begin{proof}
It is obvious that $ht(I)\leq n-1$ and the equality holds when $m=n$
and $|S_i|=1$ for all $i$, i.e., $I$ is square-free Veronese. On the
other hand, since $|S_1|+|S_2|=2(n-ht(I))\leq \sum_{i=1}^m|S_i|=n$,
we obtain that $n\leq 2\,ht(I)$. This equality holds when $m=2$ and
in this case $I=I_1*I_2$ such that each $I_i$ is generated by
$\frac{n}{2}$ distinct variables.
\end{proof}

Here, we connect matroidal ideals with graphs. Observe that if
$I\subseteq K[x_1, \dots, x_n]$ is a square-free monomial ideal of
degree $2$ then $I$ defines a simple graph $G$ with vertex set
$\{x_1, \dots, x_n\}$ and edge set $\{x_ix_j~|~x_ix_j\in I\}$. If
this is the case, we also say that $I$ is the defining ideal of $G$.
The following corollary is a consequence of Theorem~\ref{keyprop}.
\begin{coll}
\label{coll1} Let $I$ be a matroidal ideal of degree $2$ of a
polynomial ring $R=K[x_1, \dots, x_n]$. If $I$ is the defining ideal
of a simple graph $G$, then there are positive integers $t_1, \dots,
t_m$ such that $n=t_1+ \cdots +t_m$ and $G=K_{t_1, t_2, \dots,
t_m}$. In particular, if $I$ is unmixed, then $G=K_{t, t, \dots,
t}$.
\end{coll}

\begin{example}
Let $G$ be a graph defined by a matroidal ideal of degree $2$ of the
polynomial ring $R=K[x_1, \dots, x_6]$. If $G$ is unmixed, then by
Corollary~\ref{coll1}, $G=K_6$ or $K_{3,3}$ or $K_{2,2,2}$.
\end{example}



Next, we proceed to state and prove the main result in this section
which gives a characterization of unmixed matroidal ideals of degree
$d$ in a polynomial ring $R=K[x_1, \dots, x_n]$.
\begin{theo}
\label{maintheo2} Let $I$ be a matroidal ideal of degree $d$ of a
polynomial ring $R=K[x_1, \dots, x_n]$, where $n\geq 2$. If $I$ is
unmixed, then $h+d-1\leq n\leq hd$, where $h$ is the height of $I$.
In particular, $n=h+d-1$ if and only if $I$ is square-free Veronese;
and $n=hd$ if and only if $I=J_1*J_2* \cdots *J_d$, where each $J_i$
is generated by $h$ distinct variables.
\end{theo}
\begin{proof}
Observe first that the assertion holds if $d=1$ for that in this
case $I=\langle x_1\cdots x_n\rangle$ and $n=h$. Therefore we assume
that $d\geq 2$. We proceed the proof by induction on $d$. If $d=2$,
then it is the content of Corollary~\ref{d=2}. Thus we assume now
that $d\geq 3$. For $i=1, \dots, n$, let $S_i=\{\frac{u}{x_i}~|~u\in
G(I),~and~x_i\mid u\}$ and $I_i$ be the ideal generated by $S_i$.
Then $I_i$ is a matroidal ideal of degree $d-1$ with
$supp(I_i)\subseteq \{x_1, \dots, \hat{x_i}, \dots, x_n\}$ and
$$I=\sum_{i=1}^n x_iI_i.$$ We will show that $I_i$ is unmixed in
the following. For this, we prove $I_1$ for example. Let $P_1,
\dots, P_r$ be the minimal primes of $I$ that contain $x_1$ and
$Q_1, \dots, Q_s$ be the minimal primes of $I$ that do not contain
$x_1$; then $$I=P_1\cap \cdots \cap P_r\cap Q_1\cap \cdots \cap
Q_s$$ is a minimal primary decomposition of $I$. Therefore
$$x_1I_1\subseteq \langle x_1\rangle\cap I \subseteq \langle x_1\rangle\cap Q_1\cap \cdots \cap
Q_s = \langle x_1\rangle\cdot  (Q_1\cap \cdots \cap Q_s)\subseteq
x_1I_1$$ as $\langle x_1\rangle\cdot  (Q_1\cap \cdots \cap
Q_s)\subseteq I$. Hence, $I_1=Q_1\cap \cdots \cap Q_s$ which is
unmixed as $ht(Q_j)=h$ for all $j$. \par To obtain the inequality
$n\geq h+d-1$, let $t_i=|supp(I_i)|$ for $i=1, \dots, n$; then
$t_i\leq n-1$. Since $I_i$ is an unmixed matroidal ideal of degree
$d-1$ and of height $h$, $t_i\geq h+(d-1)-1$ by induction . So
$n\geq h+d-1$ as $t_i\leq n-1$. It is clear that if $I$ is
square-free Veronese then $n=h+d-1$. Conversely, if $n=h+d-1$, then
$t_i=n-1=h+(d-1)-1$, so that $I_i$ is square-free Veronese by
induction, it follows that $I$ is square-free Veronese as
$I=\sum_{i=1}^n x_iI_i$. \par To obtain $n\leq hd$, let
$T=\{x_i~|~i\neq 1, x_1x_i\mid u,~ \text{for~some}~u\in G(I)
\}\subseteq supp(I_1)\subseteq \{x_2, \dots, x_n\}$. For $i=1,
\dots, r$, choose $f_i\in Q_1\cap \cdots \cap Q_s-P_i$. Let $y\in
\{x_2, \dots, x_n\}-T$; then $x_1f_i\in I$, so that $yf_i\in
I\subseteq P_i$ by Lemma~\ref{keylem3}, it follows that $y\in P_i$
for every $i$. Therefore $h=ht(P_i)\geq 1+n-1-|T|=n-t$, where
$t=|T|$. Now, $I_1$ is an unmixed matroidal ideal of degree $d-1$
and of height $h$. By induction $h(d-1)\geq |supp(I_1)|\geq t\geq
n-h$. Therefore, $hd\geq n$. It is clear that if $I=J_1*J_2* \cdots
*J_d$ such that each $J_i$ is generated by $h$ distinct variables,
then $n=hd$. Conversely, if $n=hd$, then $P_i$ is generated by the
set $\{x_1, \dots, x_n\}-T$, so that $r=1$ and $I=P_1*(Q_1\cap
\cdots \cap Q_s)$. The assertion follows as by induction.
\end{proof}


\section{Arithmetical rank of a matroidal ideal}
The goal of this section is to study the {\it arithmetical rank} of
a matroidal ideal. For this we recall the definition of {\it
arithmetical rank} as follows.\par  Let $R$ be a Noetherian ring and
$I$ be an ideal of $R$. We say that the elements $x_1, \dots, x_m\in
R$ {\it generate I up to radical} if $\sqrt{(x_1, \dots,
x_m)}=\sqrt{I}$. The minimal number $m$ with this property is called
the {\it arithmetical rank of $I$}, denoted by $ara~ I$. If $\mu(I)$
is the minimal number of generators for $I$ and $ht(I)$ is the
height of $I$, then it is known that
$$ht(I)\leq ara~I\leq \mu(I).$$ $I$ is called {\it set-theoretic
intersection} if $ht(I)=ara~I$. The following results will be used
later.

\begin{lemma}
\label{svlem}\cite{sv} Let $P$ be a finite subset of a ring $R$.
Let $P_0, \dots, P_r$ be subsets of $P$ such that
\begin{description} \item{(i)} $\bigcup_{i=0}^r P_i=P$; \item{(ii)}
$P_0$ has exactly one element; \item{(iii)} if $p$ and $p'$ are
different elements of $P_i$ ($0<i\leq r$) there is an integer $i'$
with $0\leq i'<i$ and an element in $P_{i'}$ which divides $pp'$.
\end{description} If $q_i=\sum_{p\in P_i} p$, then
$$\sqrt{P}=\sqrt{(q_0, \dots, q_r)}.$$
\end{lemma}

\begin{lemma}
\label{*lem} Let $I$ and $J$ be two monomial ideals of $R=K[x_1,
\dots, x_n]$ such that $supp(I)\cap supp(J)=\emptyset$. Suppose
that $ara~I=u$ and $ara~J=v$. Then $ara(I*J)=u+v-1$.
\end{lemma}
\begin{proof}
See \cite[Theorem~1]{sv}, for example.
\end{proof}

\begin{lemma}
\label{qlem} Let $I$ be a matroidal ideal of degree $d$ of a
polynomial ring $R=K[x_1, \dots, x_n]$. If $ara~I\leq q(I)+1$,
then $ara~I=q(I)+1$.
\end{lemma}
\begin{proof}
From \cite{l}, we know that for $I$, the following holds:
$$pd_R~R/I\leq ara~I.$$ However, by \cite[Corollary~1.6]{ht}
$$pd_R~R/I=q(I)+1.$$ Thus, the
assertion holds.
\end{proof}

\begin{theo}
\label{maintheo4} Let $I$ be a matroidal ideal of degree $d$ of a
polynomial ring $R=K[x_1, \dots, x_n]$ and supp$(I)=\{x_1, \dots,
x_n \}$. Then $ara~I =q(I)+1$ if one of the following holds:
\begin{description}
\item{(i)} $I$ is square-free Veronese;  \item{(ii)} $I=J_1 J_2  \cdots J_d$
such that each $J_i$ is generated by $h$ distinct variables;
\item{(iii)} $d=2$,
\end{description} where $h$ is the height of $I$.
\end{theo}

\begin{proof}
By Lemma~\ref{qlem} and Theorem~\ref{maintheo3}, it is enough to
show that $ara~I\leq n-d+1$.\\ (i) Suppose that $I$ is square-free
Veronese. Let $S_i$ be the set of all square-free monomials of
degree $d$ in variables $x_1, x_2, \dots, x_i$; then $|S_d|=1$ and
$S_d\subset S_{d+1}\subset \cdots \subset S_n$. Let $P_0=S_d$ and
$P_i=S_{d+i}-S_{d+i-1}$ for $i=1, \dots n-d$; then it is easy to
check that $P=\bigcup_{i=0}^{n-d} P_i$ satisfies the assumptions of
Lemma~\ref{svlem} and $P=G(I)$. Thus, $ara~I\leq n-d+1$. \\ (ii)
Since
$ara~J_i=h$, $ara~I\leq hd-d+1=n-d+1$ follows by Lemma~\ref{*lem}.\\
(iii) By Theorem~\ref{keyprop}, we can divide the set $\{x_1,
\dots, x_n\}$ into subsets: \\
$\{x_1, \dots, x_{t_1}\},
\{x_{t_1+1}, \dots, x_{t_1+t_2}\}, \dots, \{x_{t_1+\cdots
+t_{m-1}+1}, \dots, x_n\}$ such that $n=t_1+ \cdots +t_m$ and
$x_ix_j\in I$ if and only if $i\leq t_1+\cdots +t_l<j$ for some
positive integer $l$. We may further assume that $t_m\leq
t_{m-1}\leq \cdots \leq t_1$ and arrange the generators of $I$ as
follows:
 $$ \ba{llllllll}
x_1x_{t_1+1} & \cdots &  \cdots & \cdots &
\cdots  &  \cdots & x_1x_n\\

 \cdots &  \cdots & \cdots & \cdots &   \cdots &  \cdots &   \cdots
 \\ x_{t_1}x_{t_1+1} & \cdots & \cdots  & \cdots &
\cdots & \cdots   &   x_{t_1}x_n\\
x_{t_1+1}x_{t_1+t_2+1} & \cdots & \cdots  &
\cdots &  x_{t_1+1}x_n & &   \\
\cdots & \cdots & \cdots &
\cdots &   \cdots &  &   \\
x_{t_1+t_2}x_{t_1+t_2+1} & \cdots & \cdots  &
\cdots &   x_{t_1+t_2}x_n & &   \\
\cdots & \cdots &
\cdots & \cdots &   &  &    \\
\cdots &  \cdots & \cdots & \cdots &  &  &  \\
x_{n-t_m-t_{m-1}+1}x_{n-t_m+1} & \cdots & x_{n-t_m-t_{m-1}+1}x_n &
& & &
\\ \cdots &
\cdots & \cdots &  &  &  &   \\
x_{n-t_m}x_{n-t_m+1} & \cdots & x_{n-t_m}x_n & & & &  \ea $$

From the above figure we can construct an $(t_1+\cdots
+t_{m-1})\times (t_2+ \cdots +t_m)$ matrix $A=[y_{ij}]$ with entries
in $I$ as follows: For every positive integer $i\leq t_1+\cdots
+t_{m-1}$, there is an unique nonnegative integer $k\leq m-2$ such
that $t_1+\cdots +t_k+1\leq i\leq  t_1+\cdots +t_{k+1}$. Then
$y_{ij}=x_ix_{t_1+\cdots +t_{k+1}+j}$ if $1\leq j\leq t_{k+2}+\cdots
+t_m$ and $y_{ij}=0$ otherwise. Observe that $A$ has the following
properties: \\ (a) If $y_{ij}\in  G(I)$, then $y_{ii'}\in G(I)$
whenever $i'\leq j$. \\ (b) $y_{ij}=0$, then $y_{ii'}=0$ whenever
$i'\geq j$.\\ (c) $y_{ij}=0$ whenever $i+j\geq n+1$.\\ (d) Every
generator of $G(I)$ is an entry of $A$.\\ Now let
$P_0=\{x_1x_{t_1+1}\}$ and $P_1=\{x_1x_{t_1+2}, x_2x_{t_1+1} \}$. In
general, for $0\leq l<\infty$, let
$$P_l=\{y_{ij}\in G(I)~|~i+j=l+2\}. $$ Then by (c), $P_l=\emptyset$ if $l\geq
n-1$.  Therefore by (d), $G(I)=\bigcup_{l=0}^{\infty}
P_l=\bigcup_{l=0}^{n-2} P_l$ and $|P_0|=1$. Thus, it remains to
check that the assumption (iii) of Lemma~\ref{svlem} holds. To see
this, let $y_{ij}, y_{i'j'}\in P_l$ for some $l\geq
1$. We may assume that $i<i'$. To finish the proof, we need to discuss the following two cases: \\
Case~1. $x_i$ and $x_{i'}$ are independent, i.e., $x_ix_{i'}\notin
G(I)$. In this case, let $k$ be the integer such that $t_1+\cdots
+t_k+1\leq i<i' \leq  t_1+\cdots +t_{k+1}$; then
$y_{ij}=x_ix_{t_1+\cdots +t_{k+1}+l+2-i}$ and
$y_{i'j'}=x_{i'}x_{t_1+\cdots +t_{k+1}+l+2-i'}$. Since
$l+2-i'<l+2-i$, we see that $j'<j$, it follows by (a) that
$y_{ij'}\in G(I)$. Moreover, $y_{ij'}\in P_{l'}$ for some $l'<l$
and $y_{ij'}$ divides $y_{ij}\cdot y_{i'j'}$, the assertion follows. \\
Case~2. $x_i$ and $x_{i'}$ are dependent, i.e., $x_ix_{i'}\in
G(I)$. In this case, there are two integers $k<k'$ such that
$t_1+\cdots +t_k+1\leq i\leq t_1+\cdots +t_{k+1}$ and $t_1+\cdots
+t_{k'}+1\leq i' \leq t_1+\cdots +t_{k'+1}$. Since $t_1+\cdots
+t_{k+1}<i'<n$, $1\leq i'-(t_1+\cdots +t_{k+1})\leq t_{k+2}+\cdots
+t_m$. Thus
$$\ba{ccl} x_ix_{i'} & = & x_ix_{t_1+\cdots +t_{k+1}+i'-(t_1+\cdots +t_{k+1})} \\ & = & y_{i,i'-(t_1+\cdots
+t_{k+1})} \\ &  \in &  P_{l'}\ea ,$$ where $l'=i+i'-(t_1+\cdots
+t_{k+1})-2$. Observe that that $t_1+\cdots t_{k+1}\geq i$ and
$i'+j'=l+2$ implies $i'-2<l$. We get $l'<l$. Since $x_ix_{i'}$
divides $y_{ij}\cdot y_{i'j'}$, the assertion follows.
\end{proof}

\smallskip
The following corollary is a direct consequence of
Proposition~\ref{qI} and the above theorem.
\begin{coll}\label{ci}
Let $I\subseteq K[x_1,\dots,x_n]$ be a matroidal ideal such that
supp$(I)=\{x_1, \dots, x_n \}$. Then $I$ is Cohen-Macaulay if and
only if it is a set-theoretic complete intersection.
\end{coll}

In view of Theorem~\ref{maintheo4}, we propose the following
conjecture.

\medskip \noindent{\bf Conjecture:} Let $I$ be a matroidal ideal of
degree $d$ of a polynomial ring $R=K[x_1, \dots, x_n]$. Then $ara~I
=n-d+1$.

\end{document}